\documentclass[12pt]{article}
\usepackage{amsfonts,amssymb,amsmath,eucal}
\usepackage[dvips]{graphics}
\numberwithin{equation}{section}

\newcommand{\Z}{\mathbb{Z}}
\newcommand{\C}{\mathbb{C}}

\newcommand{\T}{\Theta}
\newcommand{\la}{\lambda}

\newcommand{\vr}{\varrho}
\newcommand{\sq}{\square}
\newcommand{\sh}{{\textstyle \frac12}}

\newcommand{\lan}{\langle}
\newcommand{\ran}{\rangle}

\newcommand{\vac}{v_\emptyset}
\newcommand{\nr}[1]{:\!#1\!:}
\DeclareMathOperator{\tr}{tr}
\DeclareMathOperator{\Fr}{Fr}
\DeclareMathOperator{\sgn}{sgn}
\newcommand{\LV}{\Lambda^{\frac\infty2}V}
\newcommand{\FF}{\mathsf{F}}
\newcommand{\nF}{\nr{\FF}}
\newcommand{\TT}{\mathsf{T}}
\newcommand{\wT}{\widetilde{\TT}}
\newcommand{\M}{\mathcal{M}}
\newcommand{\KK}{\mathsf{K}}

\newcommand{\fM}{\mathfrak{M}}
\newcommand{\fT}{\sigma_{\!S}}

\newcommand{\slt}{\mathfrak{sl}_2}

\newcommand{\bla}{\boldsymbol{\la}} 
\newcommand{\fS}{\mathfrak{S}}
\newcommand{\ul}{\underline}
\newcommand{\al}{\alpha} 

\newtheorem{Theorem}{Theorem}
\newtheorem{Proposition}{Proposition}

\newtheorem{Corollary}{Corollary}

\begin{document}

\title{Infinite wedge and random partitions}
\author{Andrei Okounkov}
\date{}

\maketitle

\section{Introduction}

The aim of this paper is to show that random partitions have 
a very natural and direct connection to various structures which are
well known in integrable systems.  This
connection is arguably even more natural than, for example, in the case
of random matrices. In a sense,
we show that solitaire (which is related to increasing subsequences in
random permutations and thus to the Plancherel measure on partitions \cite{AD,BDJ})
and soliton have much more in common than the general notion of solitude. 

The other character in the title, the infinite wedge space,
is our main technical tool. Starting from the fundamental work
of the Kyoto school, this object
 plays a prominent role in integrable systems.

\subsection{}

More concretely, we consider two types of measures on partitions. The first
one is the Schur measure $\fM$ for which $\fM(\la)$ is proportional to $s_\la(x)\, s_\la(y)$,
where $s_\la$ is the Schur function and $x$ and $y$ are two independent sets 
of parameters. This measure has a natural representation-theoretic interpretation,
see Section \ref{ha}. The Plancherel measure is a very particular specialization
of $\fM$. 

The other measure is the uniform measure which gives equal weight to all
partitions of a given number. Our main findings are the following. 

\subsection{Schur measure}
We consider the correlation functions of $\fM$, that is, the probability that
the set $\fS(\la)=\{\la_i - i+\frac12\}$ contains a given set $X\subset \Z+\frac12$.
Here the $\frac12$'s are just to make certain formulas more symmetric. We prove, 
see Theorem \ref{t1},  that
\begin{equation}\label{e01}
\fM\big(\la, X\subset\fS(\la)\big) = \det \big[ K(x_i,x_j) \big]_{x_i,x_j\in K} \,,
\end{equation}
where the kernel $K$ has a nice generating function, and thus a contour integral
representation, in terms of the parameters of the Schur measure $\fM$, see Theorem \ref{t2}. 
Note the
similarity to situation in random matrices \cite{Me}, but also observe that $K$ 
does not involve any objects even remotely as complicated as polynomials 
orthogonal with respect to an arbitrary measure. 

We also prove that, as functions of the parameters of the measure $\fM$, the
correlation functions satisfy an infinite hierarchy of PDE's, namely the 
Toda lattice hierarchy of Ueno and Takasaki \cite{UT}, see Theorem \ref{TDL}. 
Again, this is a well
known phenomenon in mathematical physics that various correlation functions tend to
be $\tau$-functions of integrable hierarchies. 

Both of these results are quite straightforward once one interprets the
correlation functions as certain matrix elements in the infinite wedge
space. 

\subsection{Uniform measure}

We give a new, more simple, and conceptual proof of the formula from \cite{BO}, 
which is reproduced in Theorem \ref{t3} below,
for the following averages, called the $n$-point functions, 
$$
F(t_1,\dots,t_n) = \sum_{\la} q^{|\la|} \,  
\prod_{k=1}^n \sum_{i=1}^\infty t_k^{\la_i-i+\frac12}\,. 
$$
These $n$-point functions are sums of 
determinants involving genus $1$ theta functions
and their derivatives. The key to these $n$-point functions is
a system of $q$-difference equations \cite{BO} the representation-theoretic
derivation of which is given in Theorem \ref{t5}. 

We also investigate the local structure of a typical large partition
and find that, in contrast to the Plancherel measure case \cite{BOO},
it is quite trivial. Locally, a typical partition is a trajectory
of a random walk, with probability to make a vertical (resp.\ horizontal)
step depending on the global position on the limit shape.

\subsection{Connections and applications}

\subsubsection{Asymptotic problems}
Our formula for the correlation functions of the Schur measure 
generalizes the exact formulas
for the correlations functions for the Plancherel measure used in \cite{BOO,J}
to prove the conjecture of Baik, Deift, and Johansson \cite{BDJ} about increasing
subsequences in a random permutation and also used in \cite{BOO} to analyze
the local structure of a Plancherel typical partition in the ``bulk'' of the 
limit shape. It also generalizes
the more general formula of Borodin and Olshanski \cite{BO2} for the correlation
functions of the so called $z$-measures, see Section \ref{P&z}. The
asymptotics of $z$-measures is important for the harmonic analysis on
the infinite symmetric group $S(\infty)$, see \cite{BO1,KOV}.  Since our
kernel $K$ has a simple integral representation, the formula \eqref{e01} is
particularly suitable for asymptotic investigations. We also point out 
that our proofs are considerably simpler and arguably much more
conceptual than the ones given in \cite{BOO,BO2}.  

\subsubsection{Toeplitz determinants and Fredholm determinants}
The results of the present paper were used in \cite{BorO}
to solve a problem proposed, independently, by A.~Its and P.~Deift.  
The problem was to find a general identity of the form 
\begin{equation}\label{e02}
\det(\phi_{i-j})_{1\le i,j\le n} = \det (1-K_\phi)_{\ell_2(\{n,n+1,\dots\})}\,, 
\end{equation}
where the kernel $K_\phi$ in the Fredholm determinant admits
an integral representation in terms of the generating function $\phi(z)= \sum \phi_n z^n$
for the entries of the Toeplitz determinant.  See also \cite{BW} for alternative proofs
of the result of \cite{BorO}. Again, \eqref{e02} seems to
be particularly suitable for asymptotic analysis. 

\subsubsection{Ramified coverings, moduli spaces, and ergodic theory}
It is known, see e.g.\ \cite{D}, that the Schur and uniform measures
on partitions are related to the enumeration of ramified coverings
of the sphere and the torus, respectively. 

In the case of the torus, the exact formula \eqref{e22} was used in \cite{EO}
to compute the asymptotics of the number of ramified coverings with
given ramification type as the degree goes to infinity. These numbers,
which are certain rather complicated polynomials in Bernoulli numbers,  
can be identified with volumes of certain moduli spaces and are
important in ergodic theory \cite{EM}. For the coverings of the sphere, see \cite{O2}.

\subsection{Acknowledgments}

I very much benefited from the  discussions with S.~Bloch, A.~Borodin,
P.~Deift, A.~Eskin, M.~Kashiwara, S.~Kerov, G.~Olshanski,
Ya.~Pugaj, T.~Spencer, A.~Vershik and others. Much of this paper has grown
out of my attempt to better understand the work Borodin, Olshanski, and Kerov
on the $z$-measures. I would also
like to thank NSF for financial support
under grant DMS-9801466.

\section{Correlation functions of Schur measure}\label{s2}

\subsection{The Schur measure} 

\subsubsection{Definition} 

Consider the set of all partitions $\la$ and introduce
the following
function of $\la$ 
$$
\fM(\la)=\frac1Z\, s_\la(x)\, s_\la(y) \,
$$
where $s_\la$ are the Schur functions \cite{M}
in auxiliary variables $x_1,x_2,\dots$ and $y_1,y_2,\dots$,
and $Z$ is the sum in the Cauchy identity for the Schur
functions 
$$
Z= \sum_\la s_\la(x)\, s_\la(y) = \prod_{i,j} (1-x_i y_j)^{-1} \,.
$$
It is clear that if, for exmaple, $\{y_i\}=\overline{\{x_i\}}\subset\C$ and $Z<\infty$ then
$\fM$ is a probability measure on the set of all partitions
which we shall call the \emph{Schur measure}. This measure depends
on countably many parameters. Since our
treatment will be purely algebraic, we shall not require 
positivity from the measure $\fM$. 

\subsubsection{Harmonic analysis interpretation}\label{ha}

It is well known that the Schur functions $s_\la$ are characters of
irreducible representations of the general linear group $GL$ and also
encode characters of the symmetric group $S(n)$ as in \eqref{e007} below.
Hence one can think of $\fM$ as of a representation-valued measure such that 
$$
\fM(\la)\propto\bla\boxtimes\bla \,.
$$
where $\bla$ is a representation of $GL$ or $S(n)$ corresponding to the 
partition $\la$ and $\boxtimes$ denotes the outer tensor product.

Recall that the classical spaces of noncommutative harmonic analysis 
\begin{equation}\label{dual}
\C[\textup{Mat}(n,m)]\,,\quad  \left(\C^n\right)^{\otimes m}\!, \quad
\C[S(n)]
\end{equation}
all decompose into a multiplicity-free direct sum of representations
of the form $\bla\boxtimes\bla$ with respect to the natural actions of $GL(n)\times GL(m)$,
$GL(n)\times S(m)$, and $S(n)\times S(n)$, respectively. In either of these cases, 
$\fM(\la)$ can be interpreted as the portion of the space which transforms according
to the symmetry type $\la$. 

\subsubsection{Power-sum parameters}  

It is convenient to introduce another
parameters (also known as the Miwa variables) for the Schur measure
\begin{equation}\label{miwa}
t_k = \frac1k\sum_i x_i^k\,, \quad t'_k=\frac1k\sum_i y_i^k\,,
\quad k=1,2,\dots \,.
\end{equation}
Since the power-sum symmetric functions are free
commutative generators of the algebra of the symmetric
functions, the Schur function $s_\la(x)$ is a homogeneous
polynomial in the $t_k$'s of degree $|\la|$, where
$$
\deg t_k = k\,.
$$
More precisely, one has \cite{M}
\begin{equation}\label{e007}
s_\la(x)= \sum_{\rho=1^{r_1} 2^{r_2} 3^{r_3} \dots} \chi^\la_\rho\,
\prod_k \frac{t_k^{r_k}}{r_k!}\,,
\end{equation}
where $\rho=1^{r_1} 2^{r_2} 3^{r_3} \dots$ means that the
partition $\rho$ has $r_k$ parts of length $k$ and $\chi^\la_\rho$
is the character of a permutation with cycles $\rho$ in the
representation labeled by $\la$.
Also, we have
$$
Z=\exp\left(\sum_{k>0} k\, t_k\, t'_k\right) \,.
$$

\subsubsection{Plancherel measure and $z$-measures}\label{P&z} 

In particular, if we set
$$
t=t'=(\sqrt\xi,0,0,\dots)
$$
then $\fM$ specializes to the Poissonized Plancherel measure 
\cite{BDJ,BOO,J}
$$
\fM(\la)= e^{-\xi}\, \xi^{|\la|}
\left(\frac{\dim\la}{|\la|!}\right)^2\,,
$$ 
where $\dim\la$ is the dimension of the irreducible 
representation labeled  by $\la$ and $\xi$ is the 
parameter of the Poissonization. 

The Plancherel measure 
corresponds to taking just the dimensions of irreducible pieces
in the last space in \eqref{dual}. Similarly, taking dimensions
in the first space in \eqref{dual} and allowing $z=n$ and $z'=m$ 
be arbitrary complex one obtains the $z$-measure.  

More precisely, if one sets 
$$
t_k = \frac{\xi^{k/2} z}k \,, \quad t'_k = \frac{\xi^{k/2} z'}k\,,
\quad k=1,2,\dots\,,
$$
where $z$, $z'$, and $\xi$ are parameters, 
 then one obtains the $z$-measure \cite{KOV,BO1,BO2}
\begin{equation}\label{mz}
\fM(\la) = (1-\xi)^{zz'} \xi^{|\la|} \, 
s_\la(\underbrace{1,\dots,1}_{\textup{$z$ times}})\, 
 s_\la(\underbrace{1,\dots,1}_{\textup{$z'$ times}})\,.
\end{equation}
Here $s_\la(\underbrace{1,\dots,1}_{\textup{$z$ times}})$ stands for
the polynomial in $z$ which for positive integer values of $z$ 
specializes to $s_\la(x)$ evaluated at $x_1=\dots=x_z=1$ and $x_i=0$,
$i>z$. This polynomial is well known to be \cite{M}
$$
s_\la(\underbrace{1,\dots,1}_{\textup{$z$ times}}) =
\prod_{\sq\in\la} \frac{z+c(\sq)}
{h(\sq)}\,,
$$
where $c(\sq)$ is the content of a square $\sq\in\la$ and
$h(\sq)$ is its hook-length. More precisely, \eqref{mz}
is the mixture of the $z$-measures on partitions of a 
fixed number $n=|\la|$ by the negative binomial distribution
on $n=0,1,2,\dots$ with parameter $\xi$. 
The Plancherel measure is a limit case of the $z$-measures.

\subsection{Formula for correlation functions}\label{s23}

\subsubsection{Definition of the correlation functions}\label{s214}

It is convenient to introduce the following coordinates
on the set of partitions. To a partition $\la$ we
associate a subset
$$
\fS(\la)=\{\la_i-i+1/2\}\subset \Z+\sh \,.
$$
For example, $\fS(\emptyset)=\{-\frac 12,-\frac 32,-\frac 52\}$. The geometric
meaning of the set $\fS(\la)$ can be seen from the following picture where
the diagram of $\la$ is rotated $135^\circ$ and the elements of $\fS(\la)$ are
marked by $\bullet$
\begin{center}
\scalebox{0.5}{\includegraphics{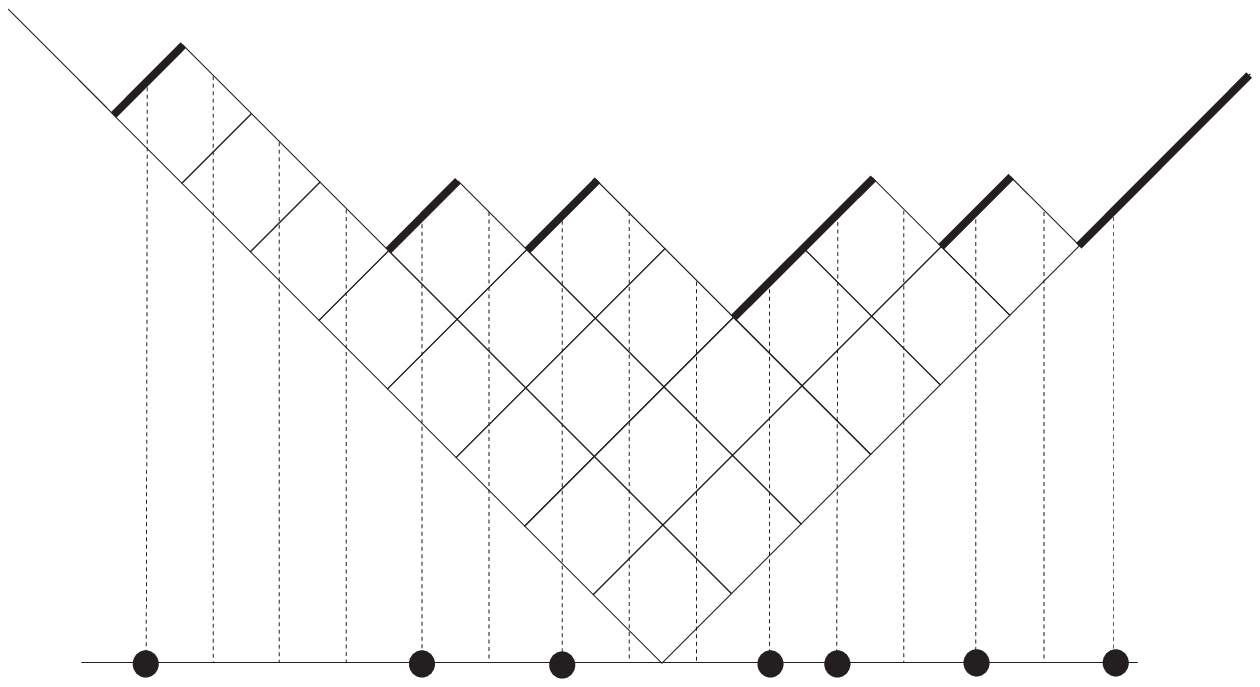}} 
\end{center}
Subsets $S\subset\Z+\frac12$ of the form $S=\fS(\la)$ can be characterized
by
$$
|S_+|=|S_-|<\infty
$$
where
$$
S_+ = S \setminus \left(\Z_{\le 0} - \sh\right) \,, \quad
S_- = \left(\Z_{\le 0} - \sh\right) \setminus S \,. 
$$
The number $|\fS_+(\la)|=|\fS_-(\la)|$ is the number of
squares in the diagonal of the diagram of $\la$ and the
finite set $\fS_+(\la)\cup \fS_-(\la) \subset \Z+\frac 12$
is known as the \emph{modified Frobenius coordinates of
$\la$}, see \cite{KO}.

Given a finite subset $X\in\Z+\sh$, define the correlation
function
$$
\rho(X)=\fM\big(\{\la, X\subset \fS(\la)\}\big) \,.
$$
Our goal in the present section  is to prove a determinantal
formula for the correlation function
$$
\rho(X)=\det\Big[K(x_i,x_j)\Big]_{x_i,x_j\in X} 
$$
where $K$ is a certain kernel which we shall
obtain momentarily. 

\subsubsection{Correlation functions as matrix coefficients}

In what follows, we assume that the reader is familiar with
the basics of the infinite wedge space and our notational conventions
as summarized for the reader's convenience in the Appendix. 

By the definition of the correlation functions $\rho(X)$ and
the equations \eqref{schur} and \eqref{e14} we have 
\begin{equation}\label{e110}
\rho(X)=\frac1Z\sum_{\fS(\la)\supset X} s_\la(x)\, s_\la(y) =  
\frac 1Z \left(\Gamma_+(t) \left( \prod_{x\in X} \psi_x \psi^*_x \right) \Gamma_-(t') \, \vac,\vac\right).
\end{equation}
We also have, see \eqref{Gamma_vac} and \eqref{Gamma_Gamma}, 
$$
\Gamma_+(t) \, \vac = \Gamma_-(t')^*\,  \vac = \vac\,, \quad  \Gamma_+(t)\, \Gamma_-(t')  =
Z\, \Gamma_-(t')  \, \Gamma_+(t) \,.
$$
{F}rom this and \eqref{e110} it follows that
\begin{equation}\label{e11}
\rho(X)=
\left( \prod_{x\in X} \Psi_x \Psi^*_x \, \vac,\vac\right)\,,
\end{equation} 
where
$$
\Psi_k = G\, \psi_k \, G^{-1} \,, \quad
\Psi^*_k = G\, \psi^*_k \, G^{-1} \,, \quad G= \Gamma_+(t)\, \Gamma_-(t')^{-1}\,.
$$
The presentation \eqref{e11} results in the following

\begin{Theorem}\label{t1} We have 
\begin{equation}\label{e12}
\rho(X)=\det\Big[K(x_i,x_j)\Big]_{x_i,x_j\in X} \,, 
\end{equation}
where $K(x,y)= \left( \Psi_x \Psi^*_y \, \vac,\vac\right)$.
\end{Theorem}

 \noindent\emph{Proof.} 
The formula \eqref{e12} can be seen as a particular case of Wick's 
theorem, see e.g.\  \cite{BS}. An equivalent direct combinatorial
argument goes as follows.  Suppose $|X|=n$ and $X=\{x_1,\dots,x_n\}$. Since
$\psi_{x_i} \psi^*_{x_j} = -  \psi^*_{x_j} \psi_{x_i}$ for $i\ne j$, 
the operators $\Psi_{x_i}$ and $\Psi^*_{x_j}$ also anti-commute. 
Therefore
$$
\rho(X)=\left(\Psi_{x_1} \cdots \Psi_{x_n}  
\Psi^*_{x_n} \cdots \Psi^*_{x_1}   \, \vac,\vac\right)\,.
$$
Applied to the vacuum $\vac$, 
the operator
$
\Psi^*_{x_n} \cdots \Psi^*_{x_1}
$
removes an $n$-tuple of the vectors $\ul{s_i}$ and then the operator
$
\Psi_{x_1} \cdots \Psi_{x_n} 
$
has to put all of them back, in $n!$ possible orders. Depending on 
parity of corresponding permutation, these $n!$ terms
appear with a $\pm$ sign and this produces the determinant. 

\subsubsection{Generating function for the kernel $K$}

Form the following generating function
\begin{align*}
\KK(z,w)&=\sum_{i,j\in\Z+\frac12} z^i w^{-j}\, K(i,j) \\
&=\left(G\, \psi(z)\, \psi^*(w) \, G^{-1} \vac,\vac\right)\,.
\end{align*}
Using the formulas \eqref{e111} we compute
$$
G\, \psi(z) \, \psi^*(w) \, G^{-1} =  
\frac{J(z)}{J(w)} \, 
 \psi(z)\, \psi^*(w) \,,
$$
where
\begin{equation}\label{J(z)}
J(z)=\frac{\gamma(z,t)}{\gamma(z^{-1},t')}=
\exp\left(\sum_{k\ge 1} t_k\, z^k - \sum_{k\ge 1} t'_k\, z^{-k}\right)=
\prod_i \frac{1-y_i/z}{1-x_i\,z } \,.
\end{equation}
{F}rom this and the formula 
$$
\left(\psi(z) \,\psi^*(w)\, \vac,\vac \right)= \sum_{j\in\Z_{\ge 0}+\frac12}
\left(\frac wz\right)^j = \frac{\sqrt{zw}}{z-w} \,, \quad |w|<|z| \,,
$$
we obtain the following

\begin{Theorem}\label{t2}  We have for $|w|<|z|$ 
$$
\KK(z,w) =  \frac{\sqrt{zw}}{z-w} \,  \frac{J(z)}{J(w)}\,,
$$
where $J(z)$ is defined by \eqref{J(z)}. 
\end{Theorem}

\noindent
The following results generalize Propositions 2.9 and
2.8 of \cite{BOO}, respectively. 

\begin{Corollary}\label{c1} Define $J_n(t,t')$ by 
$J(z)=\sum_{n\in\Z} J_n \, z^n$. We have 
$$
K(i,j)=\sum_{k\in\Z_{\ge 0}+\frac12} J_{i+k} (t,t') \, J_{-j-k} (-t,-t')  \,.
$$
\end{Corollary}

\begin{Corollary}\label{c2} We have 
$\displaystyle
\frac{\partial}{\partial t_1} \KK(z,w) = \frac{1}{\sqrt{zw}}  \frac{J(z)}{J(w)}\,. 
$
\end{Corollary}

In other words, the generating function $\KK(z,w)$ factors
after taking the partial with respect to
 $t_1$.  Similar formulas hold for partials
with respect to other parameters.

\subsubsection{$K$ as an integrable kernel}

Recall that a kernel $K$ is called integrable 
\cite{IIKS,De} if the kernel $(x-y)\, K(x,y)$
has finite rank. In more invariant terms, $K$ is integrable, if its 
commutator with the operator of multiplication by the independent
coordinate has finite rank. 

In our case, we have
\begin{multline}
\sum_{i,j\in\Z+\frac12} z^i w^{-j}\, (i-j) K(i,j) =  
\left(z\frac{\partial}{\partial z} + w\frac{\partial}{\partial w}\right) \KK(z,w)=\\
\frac{\sqrt{zw}}{z-w} \,\frac{J(z)}{J(w)} \,
\left(z\frac{\partial}{\partial z} -w\frac{\partial}{\partial w}\right) 
\log (J(z) J(w))  \,.
\end{multline}
Therefore, if the Schur measure is specialized in such a way that
$$
\frac1{z-w} \left(z\frac{\partial}{\partial z} -w\frac{\partial}{\partial w}\right) 
\log (J(z) J(w)) = \sum_{i=1}^N f_i(z)\, g_i(w)\,,
$$
for some $N$ and some functions $f_1,\dots,f_N$ and $g_1,\dots,g_N$, then
$K$ is integrable. This happens  if 
$$
\frac{\partial}{\partial z} \log J(z) \in \C(z)\,, 
$$
where $\C(z)$ stands for the field of rational functions. 

\subsection{Remarks}

\subsubsection{Multivariate Bessel functions} 
Remark that the functions $J_n(t,t')$  are the most natural multivariate generalization
of the classical Bessel functions and specialize
to the classical Bessel functions in the Plancherel specialization of 
the Schur measure \cite{BOO,J}.  The following picture, where the $z^4$ coefficient in
$\exp\left[\frac{x}{2}(z-z^{-1})+\frac{y}{2}(z^2-z^{-2})\right]$ is plotted
as a function of $x$ and $y$, illustrates how naturally yet nontrivially
this bivariate Bessel function interpolates between the classical Bessel
functions of order 2 and 4.
\begin{center}
\scalebox{0.4}{\includegraphics{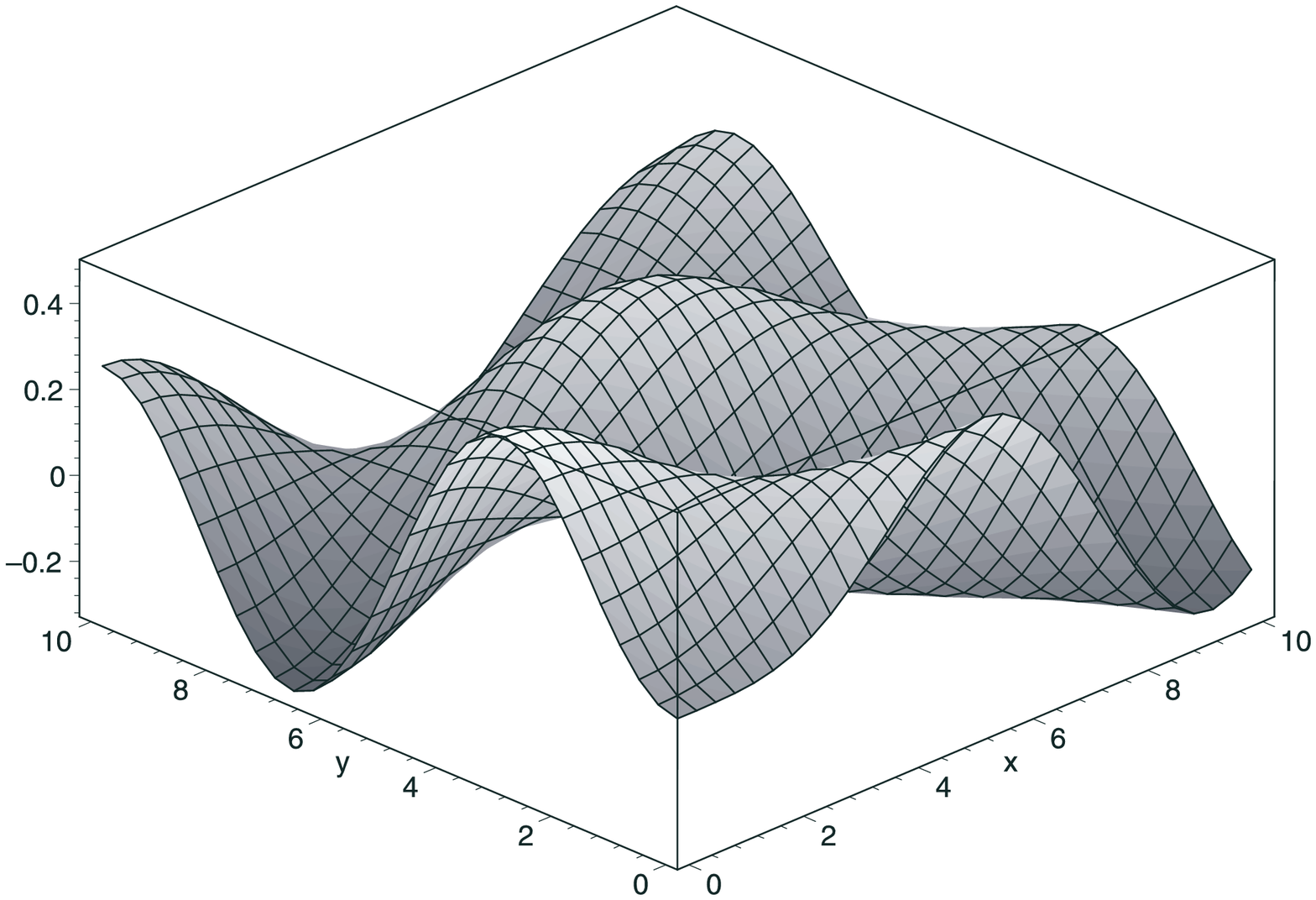}} 
\end{center}
These Bessel functions should not be confused with the (much more complicated)
multivariate Bessel functions which arise as suitable limits of the multivariate 
Jacobi polynomials. 

\subsubsection{Toda lattice hierarchy} 

The operator $\psi_k \psi_k^*$ satisfies \eqref{commOm} and so do the operators 
$$
A_X = \prod_{x\in X} \psi_x \psi_x^* \,.
$$
Denote by $X+n$ the translation of the set $X$ by an integer number $n$. We have
\begin{multline*}
(\Gamma_+(t)\, A_X  \, \Gamma_-(t) \, v_n,v_n) = 
  (\Gamma_+(t)\, R^{-n} \, A_X \, R^n \, \Gamma_-(t) \, \vac, \vac) \\ =
 (\Gamma_+(t)\, A_{X-n} \, \Gamma_-(t) \, \vac, \vac) = Z\, \rho(X-n)
\end{multline*}
Thus, we have the following 
\begin{Theorem}\label{TDL} For any set $X\subset \Z+\frac12$ the sequence
$\left\{Z\,\rho(X-n)\right\}_{n\in\Z}$ of functions of $t$ and $t'$ 
is a sequence of $\tau$-functions for the Toda lattice hierarchy.
\end{Theorem}

In the particular case when
$$
X=\{-\frac12,-\frac32,-\frac52,\dots\}\,, 
$$  
these $\tau$-functions are Fredholm determinants with the kernel $K$.
The Toda lattice equations for these Fredholm determinants are
related to the results of \cite{TW1,TW2} and of the recent paper \cite{AM}.  
See \cite{BorO} for
a discussion of the relationship between these Fredholm determinants and 
Toeplitz determinants. 

The times $t$ and $t'$ in the Toda lattice hierarchy admit the 
following combinatorial interpretation. We have
$$
\al_{-n} \, v_\mu = \sum_{\textup{$\la/\mu$ is a rim hook of size $n$}}
 (-1)^{\textup{height of $\la/\mu$ + 1}} \, v_\la \,,
$$  
where by a rim hook we mean a skew diagram $\la/\mu$ which is connected and lies
on the rim of $\la$ and the height of $\la/\mu$ is the number of rows it
occupies. The action of the operator $\Gamma_-(t')= \prod_n \exp(t'_n  \al_{-n})$ 
can be interpreted as holding our initial partition for the time $t_n$ under
a Poisson random stream of $n$-hooks, where each hook attaches to $\al$ with
a $\pm$ sign according to its height. 

Also remark that in the Plancherel case the correlations depend
only on the product $t_1 t'_1$ and so \eqref{Toda} becomes an ODE
equivalent to the usual Toda lattice \cite{T}. 

\subsubsection{Other algebras and Kerov's construction}\label{s241}
Instead of the operators $\al_n$, which generate an action of
the Heisenberg algebra, one can consider representation of 
other algebras on $\LV$. Sandwiching products of $\psi_k \psi^*_k$
between some raising and lowering operators will again
produce formulas similar to ours. The analog of the kernel $K$
will involve matrix elements of this representation. 
In particular, an $\slt$-action on $\LV$ leads again to $z$-measures
and provides another proof of the Borodin-Olshanski formula, see
\cite{O3}. The action of $\slt$ in the basis $\{v_\la\}$ gives
certain operators on partitions considered by S.~Kerov 
in his analysis of the $z$-measure \cite{Ke}.

\section{The uniform measure}\label{s3}

In this section we shall consider the uniform measure
on partitions on $n$ and the related measure
$$
\M_q(\la)= (q;q)_\infty \, q^{|\la|} \,, \quad q\in [0,1)\,, 
$$
on the set of all partitions. The normalization
factor  
$$
(q;q)_\infty=(1-q)(1-q^2)\dots
$$
makes $\M_q(\la)$ a probability measure;  this is clear
from the classical identity
$$
\sum_\la q^{|\la|} = \prod_{n\ge 1} \frac1{1-q^n} \,.
$$

\subsection{The $n$-point functions}

\subsubsection{Definition}

Again, since our treatment will be purely algebraic, we 
shall not require $\M_q$ to be positive and we shall allow 
$q$ to be any number from the unit disk $|q|<1$. 

Let $f(\lambda)$ be a function on the set of partitions. We set
$$
\langle f\rangle_q\ =(q ;q)_\infty \sum_\lambda f(\lambda)\, q^{|\lambda|}\,,
$$
assuming this expectation converges. 

Fix an integer $n\ge 1$ and variables $t_1,\ldots,t_n$. These
are not in any way related to the variables $t$ of Section \ref{s2}. 
We consider the following $n$-point function
\begin{equation}\label{e21}
F(t_1,\ldots,t_n) := \left\lan \prod_{k=1}^n\left(\sum_{i=1}^\infty
t_k^{\lambda_i-i+\frac12}\right)\right\ran_q.
\end{equation} 

This series converges if $1<|t_i|$ and $|t_1\cdots t_n|<|q|^{-1}$ and admits a 
meromorphic analytic continuation onto $\mathbb{C}^n$, see \cite{BO}
and below. The Laurent coefficients of \eqref{e21} near the point
$t_1=\dots=t_n=1$ are the averages of the form $\lan f \ran_q$ for
polynomial functions $f$ and all averages of polynomials can
be reconstructed in this fashion. 

\subsubsection{Formula for the $n$-point function}

The following formula for the $n$-point functions 
was obtained in \cite{BO}. Consider the following
genus $1$ theta function
\begin{align*}
\T(x) = \T_{11}(x;q)&=\sum_{n\in{\Z}} 
(-1)^n q^{\frac{(n+1/2)^2}2} x^{n+1/2} \\
&=q^{1/4}  (x^{1/2}-x^{-1/2}) (q;q)_\infty\, (q x;q)_\infty (q/x;q)_\infty \,. 
\end{align*}
This is the only odd genus $1$ theta function and its precise
normalization is not important because the main formula
will be homogeneous in $\T$.  Set
\begin{equation*}
\T^{(k)}(x):= \left(x\,
\frac{d}{dx}\right)^k \T(x;q)\,. 
\end{equation*} 
We have the following result from \cite{BO}

\begin{Theorem}[\cite{BO}]\label{t3} 
\begin{equation}\label{e22}
\vspace{-3 \jot}
F(t_1,\dots,t_n)= 
\sum_{\sigma\in S(n)}\, \frac 
{\displaystyle \det\left( \frac{\displaystyle \T^{(j-i+1)}(t_{\sigma(1)}\cdots
    t_{\sigma(n-j)})}{\displaystyle (j-i+1)!} 
\right)_{i,j=1}^n}
{\displaystyle \T(t_{\sigma(1)})\, \T(t_{\sigma(1)} t_{\sigma(2)}) \dots
\T(t_{\sigma(1)}\cdots t_{\sigma(n)})}
\end{equation} 
Here $\sigma$ runs through all permutations $S(n)$ of $\{1,\dotsc,n\}$, the
matrices in the numerator have size $n\times n$, and we define 
$1/(-n)!=0$ if $n\ge 1$. 
\end{Theorem}

\subsubsection{Strategy of proof} 

The proof of \eqref{e22} given in \cite{BO} is based on the following
two properties of the $n$-point functions. First, these functions 
satisfy a $q$-difference equation:
\begin{Proposition}[\cite{BO}]\label{p1} We have 
\begin{multline}\label{e23}
F(qt_1,t_2,\dots,t_n)=-q^{1/2} t_1 \dots t_n \times \\ \sum_{s=0}^{n-1} 
(-1)^s 
\sum_{1<i_1<\dots < i_s \le n} F(t_1 t_{i_1} t_{i_2} \cdots t_{i_s},
\dots,\widehat{\,t_{i_1}}, \dots,\widehat{\,t_{i_s}}, \dots) \,,
\end{multline}
where hats mean that the corresponding term should be 
omitted. 
\end{Proposition} 
Observe that since the $n$-point functions are obviously
symmetric,
we have parallel formulas for $F(t_1,\dots,qt_k,\dots,t_n)$. 

Another property of the $n$-point functions
needed in the proof of \eqref{e22} is the description of
their singularities. All singularities 
lie on the divisors
$$
q^m t_{i_1} t_{i_2} \dots t_{i_k} = 1 \,, \quad m\in\Z\,, 
$$
where $\{i_1,\dots,i_k\}\subset \{1,\dots,n\}$ is an
arbitrary subset. By symmetry, it suffices to consider
the divisor $q^m t_1 \dots t_k=1$ on which we have
\begin{equation}\label{e24} 
F(t_1,\dots,t_n)=(-1)^{m}
 \frac{q^{m^2/2} m^{k-1}}{q^m t_1\dots t_k-1}
\frac{F(t_{k+1},\dots,t_n)}{(t_{k+1}\dots t_n)^m}+ \dots \,, 
\end{equation}
where dots stand for terms regular on the divisor 
 $q^m t_1 \dots t_k=1$  and we
assume that 
$$
F(t_{k+1},\dots,t_n)=1\,,\quad k=n\,.
$$

Since we have a $q$-difference equation for the
$n$-point functions, it suffices to prove the following 
\begin{Proposition}[\cite{BO}]\label{p2} We have 
$$
F(t_1,\dots,t_n)=
 \frac{1}{t_1-1}\,
{F(t_{2},\dots,t_n)}+ \dots \,, 
$$
and $F(t_1,\dots,t_n)$ is regular on the
divisors $t_1 t_2 \cdots t_k=1$ for $k>0$. 
\end{Proposition} 

In this section  we shall give a representation-theoretic
proof of Propositions \ref{p1} and \ref{p2}. 

\subsection{Singularities of the $n$-point functions}

\subsubsection{Convergence of traces} 

 Introduce the following operators on $\LV$ 
\begin{align*}
\TT(t)&=\sum_{k\in\Z+\frac12} t^k\, \psi_k\, \psi^*_k\,, \\
\nr{\TT(t)}&=\sum_{k\in\Z+\frac12} t^k\, 
\nr{\psi_k\, \psi^*_k}\,,
\end{align*}
where the normal ordering was defined in \eqref{e112}. 
The action of the operator $\nr{\TT(t)}$ in the infinite wedge
module is well defined for any $t$. The action of $\TT(t)$
makes sense only if $|t|>1$ in which case we have
\begin{equation}\label{e25}
\nr{\TT(t)} = \TT(t) - \sum_{k=-\infty}^{-1/2} t^k  = \TT(t) - 
\frac1{t^{1/2}-t^{-1/2}} \,.
\end{equation}

Set $\fT= \sum_{k\in S_+} t^k - \sum_{k\in S_-} t^k$. It is clear that
$$
\nr{\TT(t)}\, f_S = \fT \, f_S
\,.
$$
We have the following estimates
\begin{alignat*}{2}
\left|\fT\right|& \le |t|^{\|S_+\|}+|S_-|\,,&\quad &|t|\ge 1\,,\\
\left|\fT\right|& \le |t|^{-\|S_-\|}+|S_+|\,,&\quad &|t|\le 1\,,
\end{alignat*}
where
$$
\|S\|=\sum_{k\in S} |k| \,.
$$
Since
$$
q^H \, f_S = q^{\|S_+\| + \|S_-\|} \, f_S
$$
it follows that the following trace
\begin{equation}\label{e26}
\nF(t_1,\dots,t_n)=
\tr_{\LV} \, 
\left(q^H\, \nr{\TT(t_1)}\, \nr{\TT(t_2)} \cdots  \nr{\TT(t_n)}
\right)
\end{equation}
converges absolutely provided
$$
|q| < |t_{i_1} t_{i_2} \dots t_{i_k}| < |q|^{-1} 
$$
for any subset $\{i_1,\dots,i_k\}\subset \{1,\dots,n\}$. 
Let us denote by $\Delta$ this domain of the convergence
of \eqref{e26}. 

\subsubsection{Poles of the $n$-point functions} 

Consider the trace
\begin{equation}\label{e27}
\FF(t_1,\dots,t_n)=
\tr_{\LV}
\left(q^H\, \TT(t_1)\, \TT(t_2)\cdots \TT(t_n)
\right)
\end{equation}
which converges absolutely provided $|t_i|>1$ and 
$|t_1\cdots t_n|<|q|^{-1}$. It is clear that using \eqref{e25}
one obtains a meromorphic continuation of $\FF$ onto the
domain $\Delta$. 

In particular, consider the zero charge subspace
$\Lambda_0 \subset \LV$, that is, the kernel of
the charge operator $C$. This subspace is preserved
by operators $\nr{\TT(t)}$ and $H$. It is spanned by the vectors
$v_\la$, where $\la$ is a partition, and therefore
$$
F(t_1,\dots,t_n)=(q;q)_\infty\, \tr_{\Lambda_0}
\left(q^H\, \TT(t_1)\, \TT(t_2)\cdots \TT(t_n)
\right) \,. 
$$
In fact, this is the original definition of the $n$-point
function in \cite{BO}. From the above discussion, we
have the following immediate conclusion

\begin{Theorem} The $n$-point function $F(t_1,\dots,t_n)$
admits a meromorphic continuation onto the domain $\Delta$
with simple poles along the divisors $t_i=1$ and no
other singularities. We have
$$
F(t_1,\dots,t_n)=
 \frac{1}{t^{1/2}_1-t^{-1/2}_1}\,
{F(t_{2},\dots,t_n)}+ \dots \,, 
$$
where dots stand for terms regular on the divisor $t_1=1$.
\end{Theorem}

\subsection{The $q$-difference equation}

\subsubsection{The effect of the translation operator}

Recall that $R$ denotes the translation operator \eqref{defR}. 
It is clear that 
\begin{equation}\label{e210}
R^{-1} \, \TT(t) \,R = t\, \TT(t) \,.  
\end{equation}
Therefore, from \eqref{e29} and  \eqref{decC}   we obtain 
\begin{align}
\FF(t_1,\dots,t_n)&=\sum_k \tr_{R^k \, \Lambda_0} 
\left(q^H\, \TT(t_1)\, \TT(t_2)\cdots \TT(t_n)
\right)\notag\\
&=\sum_{k} q^{\frac{k^2}2} \, (t_1 t_2 \cdots t_n)^k 
\tr_{\Lambda_0} 
\left(q^H\, \TT(t_1)\, \TT(t_2)\cdots \TT(t_n)
\right)\notag\\
&=\T_3(t_1 t_2 \cdots t_n;q)\, (q;q)_\infty \,F(t_1,\dots,t_n) \label{e211}\,,
\end{align}
where
\begin{equation}\label{e212}
\T_3(t;q)=\sum_{n\in\Z} q^{\frac{n^2}2} \, t^n \,.
\end{equation}
The function $\T_3$ satisfies the following $q$-difference
equation
$$
\T_3(qt;q)=q^{-1/2}t^{-1} \, \T_3(t;q) \,.
$$
We conclude that the $q$-difference equation \eqref{e23} for
the $n$-point function is equivalent to the following
equation
\begin{multline*}
\FF(qt_1,t_2,\dots,t_n)=\\
\sum_{s=0}^{n-1} 
(-1)^{s+1} 
\sum_{1<i_1<\dots < i_s \le n} \FF(t_1 t_{i_1} t_{i_2} \cdots t_{i_s},
\dots,\widehat{\,t_{i_1}}, \dots,\widehat{\,t_{i_s}}, \dots) \,,
\end{multline*}
which will now be established. 

\subsubsection{Proof of the $q$-difference equation}

Note that
$$
[\TT(t),\psi^*_k]=-t^k\, \psi^*_k \,.
$$
It follows that
\begin{equation}\label{e214}
\TT(t_1) \cdots \TT(t_n) \, \psi^*_k =
\sum_{P\subset\{1,\dots,n\}} (-1)^{|P|} \left(\prod_{i\in P}
t^k_i\right) \,
\psi_k^* \, \prod_{i\notin P} \TT(t_i) \,.
\end{equation}

Let us assume that $|t_i|>1$ and $|t_1\cdots t_n|<q^{-1}$ so that
the trace \eqref{e27} converges absolutely. Let us write simply
``$\tr$'' instead of $\tr_{\LV}$. From \eqref{e214} we have
\begin{multline}\label{e215}
\tr q^H \, \psi_k \TT(t_2) \cdots \TT(t_n)  \, \psi^*_k  =
\\
\sum_{P\subset\{2,\dots,n\}} (-1)^{|P|} \left(\prod_{i\in P}
t^k_i\right)\,
\tr q^H \, \psi_k\, \psi_k^* \, \prod_{i\notin P} \TT(t_i) \,.
\end{multline}
On the other hand, because of the absolute convergence, we have
\begin{align*}
\tr q^H \, \psi_k\,  \TT(t_2) \cdots \TT(t_n) \, \psi^*_k
&= \tr \psi^*_k\, q^H \, \psi_k\,  \TT(t_2) \cdots \TT(t_n) \\
&= q^k  \tr q^H \,\psi^*_k\,  \psi_k\,  \TT(t_2) \cdots \TT(t_n) \,,
\end{align*}
where we used the relation 
$
q^H\, \psi^*_k\, q^{-H} = q^{-k} \, \psi^*_k 
$\,.
Summing the equations \eqref{e215} we obtain
\begin{multline*}
\tr q^H \, \wT(qt_1)\,  \TT(t_2) \cdots \TT(t_n) =\\ 
\sum_{s=0}^{n-1} 
(-1)^{s} 
\sum_{1<i_1<\dots < i_s \le n} \FF(t_1 t_{i_1} t_{i_2} \cdots t_{i_s},
\dots,\widehat{\,t_{i_1}}, \dots,\widehat{\,t_{i_s}}, \dots) \,,
\end{multline*}
where
$$
\wT(t)=\sum_k t^k \,\psi^*_k\,  \psi_k\,.
$$
From the relations
\begin{alignat*}{2}
\TT(t) &= \nr{\TT(t)}+\frac1{t^{1/2}-t^{-1/2}}\,,& \quad &|t|>1 \,, \\
\wT(t) &= -\nr{\TT(t)}-\frac1{t^{1/2}-t^{-1/2}}\,,& \quad &|t|<1\,,
\end{alignat*}
it follows that 
$$
\tr q^H \, \wT(qt_1)\,  \TT(t_2) \cdots \TT(t_n) =
- \FF(qt_1,t_2,\dots,t_n)\,,
$$
where the right-hand side is computed by analytic continuation.
This establishes the following

\begin{Theorem}\label{t5}
 The meromorphic function $\FF(t_1,\dots,t_n)$
satisfies the following $q$-difference equation 
\begin{multline*}
\FF(qt_1,t_2,\dots,t_n)=\\
\sum_{s=0}^{n-1} 
(-1)^{s+1} 
\sum_{1<i_1<\dots < i_s \le n} \FF(t_1 t_{i_1} t_{i_2} \cdots t_{i_s},
\dots,\widehat{\,t_{i_1}}, \dots,\widehat{\,t_{i_s}}, \dots) \,.
\end{multline*}
\end{Theorem}

In view of \eqref{e211}, this establishes the $q$-difference equation 
\eqref{e23} and this concludes the proof of Propositions \ref{p1} and \ref{p2}

\subsection{An $\M_q$-typical  partition is locally a random walk}\label{s4}

The $n$-point functions encode global properties of a measure $\M_q$.
In this subsection we turn to local properties of an $\M_q$-typical
partitions and find that their asymptotics is quite trivial. More precisely,
as $q\to 1$, the local structure of an  $\M_q$-typical partition
becomes a random walk with probability to go horizontally or vertically
 depending on the global position on the rim
of the limit shape.

\subsubsection{Correlation functions for Frobenius coordinates}

Recall that the modified Frobenius coordinates of a partition $\la$
are defined as follows
$$
\Fr\la = \fS(\la)_+ \cup \fS(\la)_- \subset \Z+\sh\,,
$$
where the set $\fS(\la)=\{\la_1-1/2,\la_2-3/2,\dots\}$ was defined in \eqref{e15}.
Given a set $X\subset \Z+\sh$ define the corresponding 
\emph{correlation function} by 
\begin{align*}
\vr(X,q)&=\M_q\left(\left\{\la, X\subset\Fr\la\right\}\right)\\
&=(q;q)_\infty\, \sum_{\Fr\la\supset X}  q^{|\la|} \,,
\end{align*}
here the summation is over all $\la$ such that $X\subset\Fr\la$. 

Recall the following triple product formula for the theta function
$\T_3(z,q)$ defined in \eqref{e212}
$$
\T_3(z,q)=(q;q)_\infty \, (q^{1/2} z; q)_\infty \, (q^{1/2} /z; q)_\infty \,.
$$
Also recall that one way to prove this triple product formula uses
the identity
\begin{equation}\label{e31}
[z^0] \,  (q^{1/2} z; q)_\infty \, (q^{1/2} /z; q) = \sum_{\la} q^{|\la|}\,, 
\end{equation}
where the symbol $[z^0]$ in the left-hand side means that we take
the constant term in $z$. The combinatorial meaning of \eqref{e31} is the
following: on the right, we have all partitions, on the left, we have
all possible modified Frobenius coordinates.  The same argument establishes
the following formula. If  $X=\{x_1,x_2,\dots\}$ then 
$$
 \vr(X,q) =[z^0] \left( \prod_{x_i>0} \frac{q^{x_i}\, z}{1+q^{x_i}\, z} \, 
\prod_{x_i<0} \frac{q^{-x_i}\, z^{-1}}{1+q^{-x_i}\, z^{-1}} \,  \T_3(z,q) \right) \,.
$$
More compactly, this can be written as follows
$$
\vr(X,q) =[z^0] \left( \prod_{x\in X} \frac{(q^x\, z)^{\sgn(x)/2}}{q^{x/2}\, z^{1/2}+q^{-x/2}\, z^{-1/2}} \, 
 \,  \T_3(z,q) \right) \,.
$$
Replacing the constant term extraction by an integral we obtain
\begin{equation}\label{e32}
\vr(X,e^{-2\pi r}) = \int_{-1/2}^{1/2} \,\,\prod_{x\in X} \frac{e^{\sgn(x)\, \pi(is-rx)}}{2\cosh(\pi(is-rx))}
\, \T_3(e^{2\pi i s}, e^{-2\pi r}) \, ds\,,
\end{equation}
where $r>0$. 

\subsubsection{Asymptotics of the correlation functions}

Let $N=N(q)$ denote the expected size of a partition with
respect to the measure $\M_q$
\begin{align*}
N(q)&=(q;q)_\infty  \sum_\la |\la|\, q^{|\la|}\\
&= - q \, \frac{d}{dq} \log (q;q)_\infty = \frac 1{24} + G_2(q) \,,
\end{align*}
where $G_k(q)=\frac{\zeta(1-k)}2+\sum_n q^n \sum_{d|n} d^{k-1}$ is the 
Eisenstein series. This is, in fact, the simplest example of a polynomial
average with respect to $\M_q$ which, in the general case, are computed
by the formula \eqref{e22}. From the quasi-modular property of $G_2$ it
follows that
$$
N(e^{-2\pi r})  \sim \frac{\zeta(2)}{(2\pi r)^2}=\frac1{24 r^2} \,, \quad r\to +0 \,.
$$ 

Let us assume that the set $X$ changes with $q\to 1$ in such a way that
all limits
$$
\frac{x_i}{\sqrt N} \to  \xi_i \,, \quad x_i\in X\,, \quad q\to 1 \,,
$$
exist. That is, we consider the scaling of our typical partition by  
the factor of $\sqrt N$ in both directions, where $N$ is the area of the
typical diagram. We are interested in the asymptotics of the 
correlation functions \eqref{e32}. Observe that
by our assumption about the growth of $X$ we have
$$ 
x_i \, r \to  \frac{\xi_i}{2\sqrt 6} \,, \quad q=e^{-2\pi r} \,, \quad r\to 0 \,.
$$

The Jacobi imaginary transformation \cite{HTF} gives
$$
\T_3(e^{2\pi i s}, e^{-2\pi r}) = r^{-1/2} 
\sum_{n\in\Z} \exp\left(-\frac{\pi\,(s-n)^2}r\right) \,.
$$ 
Clearly, as $r\to +0$, only the $n=0$ summand in the above formula is
relevant and the asymptotics of the integral \eqref{e32} is determined
by the value of the integrand at $s=0$. Therefore, we obtain the following

\begin{Theorem}\label{t6} Suppose that, as $q\to1$, a finite set $X\subset\Z+\sh$
varies in such a way that all limits
$$
\frac{x_i}{\sqrt N} \to  \xi_i \,, \quad x_i\in X\,, \quad q\to 1 \,,
$$
exist, where $N=\frac 1{24} + G_2(q)$ is the expectation of $|\la|$. Then 
\begin{equation}\label{e33}
\lim_{q\to 1} \vr(X,q)  = \prod_{x_i\in X} 
\left(1+\exp\left(\frac{\pi|\xi_i|}{\sqrt{6}}\right)\right)^{-1} \,.
\end{equation}
\end{Theorem}

The error term in this formula can be also found by the Laplace
method.
This theorem says that locally a typical random
partition is just a trajectory of a  random walk. 

Theorem \ref{t6} is in agreement with Vershik's theorem \cite{V} about
the limit shape of a typical partition with respect to the uniform
measure. Namely, Vershik's theorem asserts that after the
scaling by the square root of the area in both direction, a typical
partition converges to the following limit shape
$$
\exp\left(-\frac{\pi x}{\sqrt{6}}\right)+\exp\left(-\frac{\pi y}{\sqrt{6}}\right)=1
$$
In new coordinates, $u=x-y$, $v=x+y$ this limit shape becomes
$$
v=\Upsilon(u) \,,  \quad \Upsilon(u)=
\frac{2\sqrt{6}}{\pi} \log \left(2 \cosh \left(\frac{\pi u}{2\sqrt 6}\right)\right) \,.
$$
From this we obtain the formula 
$$
\frac{1-|\Upsilon'(u)|}{2} = \left(1+\exp\left(\frac{\pi|u|}{\sqrt{6}}\right)\right)^{-1}\,,
$$
which is precisely what we see in the right-hand side of \eqref{e33}.
See the discussion in in Section 1.3 of  \cite{BOO} for why this should be the case.

\appendix

\section{Infinite wedge: summary of formulas} \label{s22}

In this section we collected, for the reader's convenience, some basic
formulas related to the infinite wedge space. This material is standard
and Chapter 14 of the book \cite{K} can be recommended as a reference.
With a few exceptions, we are closely following \cite{K}.

\subsubsection*{Definition}

Let $V$ be a linear space with basis $\left\{\ul{k}\right\}$, $k\in\Z+\sh$.
The linear space $\LV$  is, by definition, spanned by vectors 
$$
v_S=\ul{s_1} \wedge \ul{s_2} \wedge  \ul{s_3} \wedge  \dots\,,
$$
where $S=\{s_1>s_2>\dots\}\subset \Z+\sh$ is such a subset that
both sets
$$
S_+ = S \setminus \left(\Z_{\le 0} - \sh\right) \,, \quad
S_- = \left(\Z_{\le 0} - \sh\right) \setminus S 
$$
are finite. We equip $\LV$ with the inner product 
in which the basis $\{v_S\}$ is orthonormal. 

\subsubsection*{Free fermions}

Introduce the following operators in $\LV$. 
The operator $\psi_k$ is the exterior multiplication by $\ul{k}$
$$
\psi_k \left(f\right) = \ul{k} \wedge f  \,. 
$$
The operator $\psi^*_k$ is the adjoint operator.
These operators satisfy the canonical anti-commutation relations
$$
\psi_k \psi^*_k + \psi^*_k \psi_k = 1\,,
$$
all other anticommutators being equal to $0$. It is clear that
\begin{equation}\label{e14}
\psi_k \psi^*_k \,\, v_S = 
\begin{cases}
v_S\,, & k \in S \,, \\
0 \,, & k \notin S \,.
\end{cases}
\end{equation}
Introduce the following generating functions
$$
\psi(z)=\sum_{i\in\Z+\frac12} z^i\, \psi_i \,, \quad
\psi^*(w)=\sum_{j\in\Z+\frac12} w^{-j}\, \psi^*_j \,.
$$ 
and the normally ordered products
\begin{equation}\label{e112}
\nr{\psi_k\, \psi^*_k} =
\begin{cases}
\psi_k\, \psi^*_k\,, & k>0 \,,\\
-\psi^*_k\, \psi_k\,, & k<0 \,,
\end{cases}
\end{equation}

\subsubsection*{Energy, charge, and translation operators} 

Define the energy and charge operators by
$$
H=\sum_k k \nr{\psi_k\, \psi^*_k}\,, \quad
C=\sum_k \nr{\psi_k\, \psi^*_k}\,.
$$
We have
\begin{equation}\label{actC}
C\, v_S = \left(|S_+|-|S_-|\right) \, v_S\,.
\end{equation}
We will call the operator $R$ defined by 
\begin{equation}\label{defR}
R \, \ul{s_1} \wedge \ul{s_2} \wedge  \ul{s_3}  \wedge  \dots =
\ul{s_1+1} \wedge \ul{s_2+1} \wedge  \ul{s_3+1}  \wedge  \dots 
\end{equation}
the translation operator. In \cite{K}, it is denoted by $q$. Clearly
$$
R\,\psi_k \, R^{-1} = \psi_{k+1} \,, 
\quad  R\,\psi^*_k \, R^{-1} = \psi^*_{k+1}\,.
$$
It follows that
\begin{align}
R^{-k} \, C \, R^k &= C+k \,, \notag\\
R^{-k} \, H \, R^k &= H+kC+\frac{k^2}2 \label{e29}\,.
\end{align}

We have the following charge decomposition 
\begin{equation}\label{decC}
\LV=\bigoplus_{k\in\Z} R^k \, \Lambda_0 \,,
\end{equation}
where $\Lambda_0=\ker C$ is the zero charge subspace.  From \eqref{actC} we see
that $\Lambda_0$ is spanned by the following vectors
\begin{equation}\label{e15}
v_\la=v_{\fS(\la)}\,, \quad \fS(\la) = \{\la_i-i+\tfrac12\} \,, 
\end{equation}
where $\la$ is a partition. The vacuum vector
$$
\vac = \ul{-\tfrac12} \wedge \ul{-\tfrac32} \wedge \ul{-\tfrac52} \wedge \dots
$$
is  the vector with the minimal 
eigenvalue of $H$. The vector
\begin{equation}\label{vacn}
v_k=R^k \, \vac = \ul{k-\tfrac12} \wedge \ul{k-\tfrac32} \wedge \ul{k-\tfrac52} \wedge \dots\,,
\quad k\in\Z
\end{equation}
is the vacuum vector in the charge $k$ subspace. 

\subsubsection*{Bosons and vertex operators} 
Define the operators $\al_n$ by 
$$
\al_n = \sum_{k\in\Z+\frac12} \psi_{k-n}\, \psi^*_k \,, \quad n=\pm 1,\pm2, \dots \,. 
$$
They satisfy the Heisenberg commutation relations
$$
\left[\al_n, \al_m\right] = n \, \delta_{n,-m} \,,
$$
see the formula (14.10.1) in \cite{K}. 
Evidently, $\al^*_n=\al_{-n}$. It is clear from definitions that 
\begin{equation}\label{e107}
[\al_n,\psi(z)]=z^n\, \psi(z)\,, \quad  [\al_n,\psi^*(w)]=- w^{n}\, \psi^*(w) \,.
\end{equation}

Given any sequence $s=(s_1,s_2,\dots)$, define 
$$
\Gamma_{\pm}(s) = \exp\left(\sum_{n=1}^\infty s_n \, \al_{\pm n} \right)  \,. 
$$
Note that this is slightly different from the definition of $\Gamma_\pm$ in \cite{K}. We have 
\begin{equation}\label{Gamma_vac}
\Gamma_+(s) \, v_m = v_m  \,. 
\end{equation}
Also observe that $\Gamma_{\pm}^*=\Gamma_{\mp}$ and 
\begin{equation}\label{Gamma_Gamma}
\Gamma_+(s)\, \Gamma_-(s') = e^{\sum n\, s_n s'_n} \, \Gamma_-(s')\, \Gamma_+(s) \,.
\end{equation}
We have from \eqref{e107} 
\begin{align}
\Gamma_\pm(s) \, \psi(z)  &= \gamma(z^{\pm 1},s) \, \psi(z) \, \Gamma_\pm (s) \,, \notag  \\ 
\Gamma_\pm(s) \, \psi^*(z)  &= \gamma(z^{\pm 1},s)^{-1} \, \psi^*(z) \, \Gamma_\pm (s)\,. \label{e111}
\end{align}
Here $\gamma(z,t)$ is the following generating function
\begin{equation}\label{e109} 
\gamma(z,t)= \exp\left(\sum_{n\ge 1} t_n\, z^n\right) = \prod_i \frac1{1-x_i z} = 
\sum_{n\ge 0} z^n\, h_n(x) \,,
\end{equation}
where the variables $x$ are related to $t$ by \eqref{miwa} and the $h_n$'s are the complete
homogeneous symmetric functions.

The fermions $\psi(z)$ and $\psi^*(z)$ have the following expression 
in terms of $\Gamma_\pm$ (see Theorem 14.10 in \cite{K})
\begin{align}
\psi(z) &= z^{C} \, R\, \Gamma_-(\{z\}) \, \Gamma_+(-\{z^{-1}\}) \,, \notag \\
\psi^*(z) &= R^{-1} \, z^{-C} \Gamma_-(-\{z\}) \, \Gamma_+(\{z^{-1}\}) \,, \label{psi_Ga} 
\end{align}
where $\{z\}=(z,\frac {z^2}2,\frac {z^3}3\, \dots )$. 

\subsubsection*{Schur functions as matrix elements}

Using \eqref{e111}, one checks that
\begin{equation}
(\Gamma_-(t) \, v_\mu, v_\la) = \det(h_{\la_i-\mu_j+j-i}(x))  = s_{\la/\mu}(x) \,,
\end{equation}
where $s_{\la/\mu}$ is the skew Schur function and the second equality is the Jacobi-Trudy
identity for $s_{\la/\mu}$, see \cite{M}. In particular
\begin{equation}\label{schur}
\Gamma_-(t) \, \vac = \sum_\la s_\la(x) \, v_\la \,.
\end{equation}

\subsubsection*{Toda lattice and Hirota equations}

Let $A$ be an operator on $\LV$ such that
\begin{equation}\label{commOm}
\left[A\otimes A, \Omega\right] = 0 \,, \quad \Omega=\sum \psi_k \otimes \psi^*_k \,.
\end{equation}
Then the following functions of $t$ and $t'$ 
\begin{equation}\label{tauA}
\tau_n(t, t') = \left(\widetilde{A}\, v_n,v_n\right)\,, \quad  \widetilde{A}=\Gamma_+(t) \, A \, \Gamma_-(t') \,, 
\end{equation}
satisfy a hierarchy of Hirota equations. If the operator
$A$ satisfies \eqref{commOm}, then so does $\widetilde{A}$. This commutation equation can be 
rewritten as
\begin{equation}\label{commOm2}
[z^0] \left(\widetilde{A}\otimes \widetilde{A}\, \right) \left(\psi(z)\otimes \psi^*(z)\right) = 
[z^0] \left(\psi(z)\otimes \psi^*(z)\right) \left(\widetilde{A}\otimes \widetilde{A}\,\right) \,,
\end{equation}
where $[z^0]$ stands for the constant term in $z$. Taking the matrix coefficients
of \eqref{commOm2} between the vectors
\begin{align*}
\Gamma_-(s')\,  v_m \otimes \Gamma_-(-s')\,  v_{l+1} \,, \quad
 \Gamma_-(s)\, v_{m+1} \otimes \Gamma_-(-s)\,  v_{l}\,,
\end{align*}
and using the equations \eqref{Gamma_vac}, \eqref{Gamma_Gamma}, and
\eqref{psi_Ga}, one obtains the relation
\begin{multline}\label{Hirota}
[z^{l-m}] \, \gamma\left(\tfrac1z,-2s'\right)\, \tau_{m+1}(t+s,t'+s'+\{z\}) \, \tau_{l} (t-s,t'-s'-\{z\}) =\\
[z^{m-l}] \, \gamma\left(\tfrac1z, 2s\right)\, \tau_{m} (t+s-\{z\},t+s') \, \tau_{l+1} (t-s+\{z\},t-s') \,,
\end{multline}
which is satisfied for any $m,l\in\Z$ and any pair of sequences $s$ and $s'$. Up to the transformation
$t\to -t$, these are the Hirota equations for the Toda lattice hierarchy
of Ueno and Takasaki, see Theorem 1.11 in \cite{UT}.
Expanding \eqref{Hirota} into the Taylor series in $s$ and $s'$, one obtains a system of PDE's
satisfied by the $\tau_n$'s.  In particular, setting $l=m-1$ and 
extracting the coefficient of $s'_1$ gives
\begin{equation}\label{Toda}
\frac{\partial^2}{\partial t_1 \, \partial t'_1} \log \tau_n = 
\frac{\tau_{n+1}\, \tau_{n-1}}{\tau_n^2} \,,
\end{equation}

It may be helpful to point out that these Hirota equations are, essentially, 
Pl\"ucker-type relations for the infinite determinants which $\tau$-functions are
by their definition. In particular, \eqref{Toda} corresponds to the famous 
determinantal identity which is often associated with Lewis Carroll \cite{LC}, 
but was first established by P.~Desnanot in 1819,
see \cite{muir} (I am grateful to A.~Zelevinsky for this historical information.)

\end{document}